\documentclass[12pt]{article}
\newtheorem{thm}{Th\'eor\`eme}[section]
\newtheorem{prop}{Proposition}[section]

\newtheorem{defin}{Definition}[section]

\begin{document}

\title{\Large \bf R\'ealisation de big\`ebres
(III).\\Module induit,big\`ebre enveloppante, dual.
L'exemple des champs de vecteurs sur $C^{p}$
et leur int\'egration. }
\author{Eric Mourre .}
\maketitle

\begin{abstract}
In this article an interpretation and a proof of some
classical \\theorems in analysis on the integration
of analytic vectors fields are derived from the algebraic
method of realization of bialgebras which are constructed
with the data of a linear application from a coalgebra 
into the algebra  of right (or left) invariant operators on
an approximated coalgebra 
\cite{mo1}, \cite{mo2},\cite{mo3}. The results are obtained
from these general algebraic construction and theorems by
introducing the more restrictive notion of induced module.
Then the associated envelopping bialgebra is defined and
naturally belongs to the dual of the tensor algebra over the
approximated coalgebra .\\ An interesting technical
contribution is due to the "coproduct", or coproducts
given by the approximated coalgebra, in the classical case
of vectors fields at least. 

\end{abstract}
 Centre de Physique Theorique ,C.N.R.S. Luminy case 907,
13288 Marseille Cedex 9, France ;
U.M.R. 6207. CPT-P57-2006
\section{Introduction .}
Dans les articles pr\'ec\'edents \cite{mo1},\cite{mo3} une
m\'ethode purement alg\'ebrique  
de r\'ealisation de big\`ebres est pr\'esent\'ee dans
un cadre plut\^ot g\'en\'eral , et en suite un th\'eor\`eme
de dualit\'e est d\'emontr\'e sous des hypoth\`eses plus
restrictives, adapt\'ees \`a la r\'ealisation de
big\`ebres  telles que par exemple les g\'en\'eralisations
bien connues de l'alg\`ebre enveloppante d'une alg\`ebre de
Lie de dimension finie \cite{ Ma}.   Cet article a pour but
de montrer comment l'on peut d\'efinir, construire, et
utiliser l'alg\`ebre enveloppante des champs de vecteurs sur
$C^{p}$, et son dual, et retrouver la th\'eorie de Cauchy
pour l'int\'egration des champs de vecteurs analytiques sur
$C^{p}$ , \`a partir de la m\'ethode g\'en\'erale de
r\'ealisation de big\`ebres et du th\'eor\`eme de dualit\'e. 
En particulier on utilisera la notion de cog\`ebre
approch\'ee: $F_{a}$, dual restreint de limites inductives
particuli\`eres d'alg\`ebres de dimensions finies, ainsi que
la notion d'op\'erateurs invariants \`a droite, ou \`a
gauche, sur l'alg\`ebre tensorielle 
$T(F_{a})$ construite sur $F_a $, notions et outils  
introduits dans la version r\'evis\'ee  \cite{mo2} (plus
alg\'ebrique).\\
\\La partie th\'eorique de cet article
est restreinte \`a l'introduction du probl\`eme
int\'eressant de l'induction de relations sur un 
module alg\'ebrique de 
$T(F_a )$
 pour l'action d'une classe d'op\'erateurs invariants 
\`a gauche (paragraphe2). Cette approche g\'en\'erale 
permettra  d'interpr\'eter une classe d'op\'erateurs
invariants \`a gauche comme champs de vecteurs sur un module
alg\'ebrique que les relations induites rendent alors
hom\'eomorphe \`a l'alg\`ebre des fonctions polynomiales sur
$C^{p}$ (paragraphe 3 ); la partie technique est
d\'etaill\'ee, dans ce cadre plut\^ot g\'en\'eral , pour le
th\'eor\`eme d'existence concernant l'int\'egration des
champs de vecteurs analytiques, l'interpr\'etation, et
l'utilisation du coproduit  sur la cog\`ebre approch\'ee $F_a
$.\\ En ce qui concerne la discussion th\'eorique sur les
modules induits et la construction de modules particuliers
admettant   une notion de champs de vecteurs non triviale,
on a \'et\'e amen\'e \`a supposer des hypoth\`eses fortes
 sur la structure de leurs
relations .

\section{Cog\`ebre approch\'ee, op\'erateurs invariants
\`a droite sur son alg\`ebre tensorielle, et r\'ealisations
de big\`ebres.}
\subsection{Cog\`ebres approch\'ees adapt\'ees \`a la
r\'ealisation des  champs de vecteurs sur
$C^{p}$.} Soit $V$ l'espace vectoriel $C^{N^{p}}$, les
matrices
\'el\'ementaires sur cet espace vectoriel serons not\'ees
$E_{n}^{m}$ o\`u $n=(n_1 ,n_2 ,..,n_p)$ et
$m=(m_1,m_2,..,m_p )$, sont des multi-indices \`a p
composantes. Sur l'espace vectoriel de ces matrices on a une
structure d'alg\`ebre et on appelle $E$ la limite inductive
obtenue par les injections des sous espaces vectoriels
$C^{(0,n_1)\times (0,n_2 )\times ..\times (0,n_p)}$  dans
$C^{(0,m_1)\times (0,m_2 )\times ..\times (0,m_p)}
$ lorsque $n_{i} \leq m_{i} $ pour tout $i\in (1,2,..,p)
$.\\ Le dual restreint de cette limite inductive
particuli\`ere d'alg\`ebres est muni d'une structure de
cog\`ebre approch\'ee \cite{mo2}  ; en
particulier une cog\`ebre approch\'ee est aussi , comme
espace vectoriel, une limite inductive de cog\`ebres de
dimensions finies ; on note cette cog\`ebre approch\'ee par
$F_{a}$ et ses \'el\'ements:
$f_{n}^{m} $, n et m \'etant des multi-indices \`a p
\'el\'ements . Le coproduit \`a priori formel est d\'efini
par :
$$ \Delta_{F} f_{n}^{m}=
\sum_{k\in N^{p}}
f_{k}^{m} \otimes
f_{n}^{k} 
$$
Dans ce qui suit on notera les multi-indices \`a p
\'el\'ements , $n$ et on d\'esignera par $\vert n \vert =
\sup_{i\in (1,..,p)} 
n_{i} $ .

\begin{defin}{Op\'erateurs
invariants \`a droite ou \`a gauche r\'eguliers, sur la
cog\`ebre approch\'ee $F_a $ .}
\end{defin}
Un op\'erateur $O_{r}$ ou matrice d'\'el\'ements $O_{n}^{m}$
est dit r\'egulier si et seulement si il existe $ c\in N$
tel que 
$O_{n}^{m}=0$ quand $ \vert n-m\vert > c . $\\
Les op\'erateurs invariants \`a droite, resp. invariants \`a
gauche r\'eguliers, sont les op\'erateurs dans $End(F_a) $
obtenus par des op\'erateurs r\'eguliers de la fa\c{c}on
suivante :
$x\in Invd,r(F_a) $ :
$$x(f_{n}^{m}) = \sum_{k} f_{k}^{m} (O_r).f_{n}^{k} ,$$ 

et $x\in Invg,r(F_a) $ :
$$x(f_{n}^{m})=\sum_{k} f_{k}^{m}.f_{n}^{k}(O_r ). $$
Le fait de ne consid\'erer que des contractions du coproduit
par des op\'erateurs r\'eguliers ram\`ene les sommations \`a
priori infinies \`a des sommations finies . \\
Remarquons que la cog\`ebre approch\'ee $F_a $ est bien
munie d'une counit\'e (\`a droite et \`a gauche): $\epsilon
_{F_a}$, correspondant aux contractions du coproduit formel
par l'op\'erateur identit\'e. \\ Pour introduire la notion
de champs de vecteurs , dans un cadre g\'en\'eral, on aura
besoin d'induire des relations sur un sous module de $T(F_a)
$ par des choix particuliers; pour cela, de m\^eme que pour
\'etudier partiellement l'alg\`ebre enveloppante associ\'ee
au champs de vecteurs sur $C^p $ on aura besoin des deux
th\'eor\`emes g\'en\'eraux suivants qui sont essentiellement
expos\'es de mani\`ere alg\'ebrique dans un cadre non 
co-commutatif, dans les articles 
\cite{mo1}, \cite{mo2}, \cite{mo3} . Le th\'eor\`eme 2.1
suivant est en fait une l\'eg\`ere extension du  th\'eor\`eme
correspondant
 de \cite{mo2} parce qu'ici on
utilise les op\'erateurs invariants \`a droite construits
avec des op\'erateurs r\'eguliers. D'autre part on utilise
implicitement l'alg\`ebre tensorielle construite sur une
limite inductive : voir \cite{bou}.
\subsection{R\'ealisation de big\`ebres.} 
\begin{defin}
Une cog\`ebre $L$ est un espace vectoriel muni d'un coproduit
$\Delta_{L} $ coassociatif et d'une counit\'e 
$\epsilon_{L}$ (
\`a droite et \`a gauche ).

\end{defin}

\begin{thm}
Soit une application lin\'eaire $ x:L\rightarrow Invd,r(F_a
)$  de la cog\`ebre $L$ dans les
op\'erateurs invariants \`a droite r\'eguliers sur la
cog\`ebre approch\'ee $F_a $ .
Alors:\\
A) Il existe une unique application linaire 
$$X:L \rightarrow Invd(T(F_a)) \subset End(T(F_a)) ,$$
qui v\'erifie :\\
1) $ X(l) (1) =\epsilon_{L}(l).1 $, pour $1\in T(F_{a})$ ;\\
2) $X(l) (f) = x(l) (f) $, $\forall f\in F_a $ ;\\
3) pour tout $w_1 , w_2 \in T(F_a ) $ ,
$X(l)(w_1 .w_2 ) =\sum_{k}X(l_{k}^{'})(w_1 )
.X(l_{k}^{''})(w_2 ) ,$\\
o\`u 
$\Delta_{L}(l)= \sum_{k}l_{k}^{'}\otimes l_{k}^{''} ;$\\
4) pour tout $n\geq 0$ , $X(l): \otimes^{n}F_a \rightarrow
\otimes^{n}F_a $
\\ 5)les op\'erateurs $X(l)$ sont des op\'erateurs
invariants \`a droite sur l'alg\`ebre tensorielle $T(F_a )$.
 \\ B) L'alg\`ebre $U_{x}\subset End(T(F_a ))$ engendr\'ee
par les op\'erateurs $X(l) $ et l'identit\'e est munie d'une
unique structure de big\`ebre
$U_{x}(\Delta_{L},\epsilon_{L})$, o\`u\\ $\Delta _{L}:U_x
\rightarrow U_x \otimes U_x $ \'etend par
morphisme d'alg\`ebres le coproduit d\'efini sur les
g\'en\'erateurs: $\Delta_{L}(X(l)) =\sum_{k}
X(l_{k}^{'})\otimes X(l_{k}^{''})$ , $\Delta id = id\otimes
id$, et o\`u la counit\'e $\epsilon _{L}$ est le morphime
d'alg\`ebres v\'erifiant : $\epsilon_{L}(id)=1 ,\  
\epsilon_{L}(X(l)) =\epsilon_{L}(l) $ .

\end{thm}
Remarque1) Le th\'eor\`eme est \'ennonc\'e pour les
op\'erateurs invariants \`a droite, mais reste
valable pour les op\'erateurs invariants \`a gauche puisque
pour s'y ramener il suffit de consid\'erer la cog\`ebre
approch\'ee oppos\'ee de la cog\`ebre approch\'ee $F_a $.\\
remarque 2) Le th\'eor\`eme ci-dessus est 
g\'en\'eral, puisqu'il n'y a pas de restrictions sur la
cog\`ebre $L$, ni sur la forme particuli\`ere des
op\'erateurs invariants \`a droite r\'eguliers $x(l)$ . Pour
introduire dans ce cadre l'alg\`ebre enveloppante
associ\'ee aux champs de vecteurs sur $C^p$, la cog\`ebre
$L$ sera choisie comme \'etant une
cog\`ebre de Leibnitz, c'est \`a dire le dual d'une
alg\`ebre
$K$ poss\'edant une unit\'e $k_0 $ et des \'el\'ements $k_i$
v\'erifiant: $k_i .k_j =0 $ pour $i\  et\ j \neq 0$. De
plus pour pouvoir parler de champs de vecteurs sur
$C^{p}$ il nous faudra   induire  (gr\^ace au th\'eor\`eme de
dualit\'e, et aux caract\'erisations des relations qu'il
fournit), par le choix des op\'erateurs
$x(l_{i})
$, des relations sur un sous module commun de $ T(F_a)$
telles, qu'il soit isomorphe \`a l'alg\`ebre des fonctions
polynomiales  sur $C^{p}$ .\\ On
verra l'interpr\'etation et l'utilit\'e du coproduit sur la
cog\`ebre approch\'ee
$F_{a}$. Pour cela, ainsi que pour rester dans le cadre plus
g\'en\'erale d'alg\`ebres enveloppante non co-commutatives on
aura besoin du th\'eor\`eme de
dualit\'e dont nous allons d\'ecrire le r\^ole et donner l'
\'ennonc\'e dans le paragraphe suivant .
  
\subsection{Th\'eor\`eme de dualit\'e .}

\begin{defin}Une cog\`ebre de type fini est une cog\`ebre
telle que tout sous espace vectoriel de dimension fini soit
contenu dans une sous cog\`ebre de dimension finie.

\end{defin}
Soit $ x:L \rightarrow Invd,r(F_{a})$, une application
lin\'eaire d'une cog\`ebre de type fini dans les op\'erateurs
invariants \`a droite r\'eguliers sur la cog\`ebre
approch\'ee $F_a $, et $X(l)$ les op\'erateurs invariants
\`a droite agissant sur $T(F_{a})$, $$X(l)\in
Invd(T(F_{a}))
\subset End(T(F_{a}))\ ,$$ donn\'es par le th\'eor\`eme 2.1,
et $U_x (\Delta_{L},\epsilon_{L})$ la big\`ebre engendr\'ee
par les op\'erateurs $X(l)$ et l'identit\'e dans
$Invd(T(F_{a}))$. \\ Cette construction nous donne donc une
repr\'esentation de l'alg\`ebre $T(L)$ sur $T(F_{a})$ : un
morphime $\pi_{x}:T(L) \rightarrow Invd(T(F_{a}))$; de plus
cette repr\'esentation est une action de la big\`ebre libre
$T(L)$ sur l'alg\`ebre $T(F_{a})$ c'est \` a dire:\\
$w\in T(L) \ ,  et\  z_1 , z_2 \in T(F_{a})$ nous avons :
$$\pi_{x}(w)(z_1 .z_2 ) = \sum_{k}\pi_{x}(w_{k}^{'})(z_1
).\pi_{x}(w_{k}^{''})(z_2 )\ ,\ \ \Delta_{L}(w)=
\sum_{k}w_{k}^{'}\otimes w_{k}^{''}\ .$$
Par d\'efinition une relation de l'alg\`ebre $U_{x}$ est un
\'el\'ement $w$ de l'alg\`ebre $T(L)$ tel que
$\pi_{x}(w)=0\in U_{x}\subset End(T(F_{a})) \ .$\\
\\
Etant donn\'ee l'application lin\'eaire $ x:L \rightarrow
Invd,r(F_{a})$, d'une cog\`ebre de type fini dans les
op\'erateurs invariants \`a droite r\'eguliers sur la
cog\`ebre approch\'ee $F_a $, nous allons \'etudier gr\^ace
au th\'eor\`eme de dualit\'e la structure des relations
induites sur $T(F_a )$ par la donn\'ee de cette application
$x$ et plus pr\'ecisemment, les relations induites sur un
module alg\'ebrique.\\
\begin{defin}

 Les modules alg\'ebriques de $F_{a}$ pour
l'action des op\'erateurs invariants \`a droite  sont fournis
par les vecteurs colonnes de la cog\`ebre approch\'ee
$f_{m}^{n}$,
$m$ \'etant fix\'e; pour l'action des op\'erateurs invariants
\`a gauche les modules alg\'ebriques \'etant les
vecteurs lignes . 
\end{defin}
\vspace{0.2cm}
{\bf Remarque }
Dans cet article on se limitera dans la pratique aux seules
cog\`ebres approch\'ees d\'ecrites plus haut, c.a.d
: $F_{a}=(f_{n}^{m},\ n , m \in N^{p})$, $p$ \'etant
fix\'e; cette restriction est en fait inutile puisque pour
$p$ arbitraire les \'el\'ements font eux m\^emes partie du
dual restreint d'une double limite inductives d'alg\`ebres de
dimensions finies avec unit\'es. De plus il faut mentionner 
la stabilit\'e par somme directe et produit tensoriel 
des limites inductives particuli\`eres d'alg\`ebres
dont une restriction du dual d\'efinissent les cog\`ebres
approch\'ees; on a donc la m\^eme stabilit\'e pour les
cog\`ebres approch\'ees.\\
\\ 
Pour introduire le r\^ole du th\'eor\`eme de dualit\'e dans
le  probl\`eme de l'induction de relations int\'eressantes
sur un module alg\'ebrique d'une cog\`ebre approch\'ee 
$F_a $, consid\'erons donc la donn\'ee 
g\'en\'erale d'une application lin\'eaire\\ $x:L\rightarrow
Invg,r(F_{a})$, o\`u $L$ est une cog\`ebre de type finie
et $U_x(\Delta_{L},\epsilon_{L}) $ la big\`ebre contruite
par le th\'eor\`eme [2.1] et $\pi_x :T(L) \rightarrow U_x
\subset Invg (T(F_a )) $ le morphisme de big\`ebres alors
d\'efini.       Soit
$L_0
$, une sous cog\`ebre de dimension finie de $L$; la
restriction de
$x$
\`a la sous cog\`ebre $L_0 $, nous donne une application $
x_0 :L_0
\rightarrow Invg,r(F_a ).$\\ Par d\'efinition pour $l\in L_0
$ nous avons : 
$$x_{0}(l) (f_{n}^{m}) =\sum_{k}f_{k}^{m}.f_{n}^{k}(O_{r}(l))
\ .$$ Consid\'erons donc l'application:
$x_{0,t} :l\rightarrow O_{r}(l) $;
par transposition on obtient une application :
$y_{0} :F_a \rightarrow K_0$, o\`u $K_0 $ est l'alg\`ebre
duale de la cog\`ebre $L_0 $. \\
La restriction de cette application \`a toute cog\`ebre de
dimension finie de la forme:\\ $( f_{n}^{m} \ , \vert n\vert
 ,\
\vert m\vert \leq q )$ est bien d\'efinie . Mais pour
comprendre la nature des relations accessibles sur $T(F_{a})$
il nous faudra utiliser le th\'eor\`eme de dualit\'e,
voir \cite{mo3}, suivant qui d'ailleurs ne s'applique pas
directement m\^eme  pour l'\'etude des relations entre les
seuls
\'el\'ements $ f_{n}^{m} \ , \vert n\vert
 ,\
\vert m\vert \leq q $.\\ La version du th\'eor\`eme
de dualit\'e suivant, concernant les applications lin\'eaires
d'une cog\`ebre de dimension finie  dans les op\'erateurs
invariants
\`a gauche sur une cog\`ebre de dimension finie, reste
valable, puisque les op\'erateurs invariants \`a gauche sont
des op\'erateurs invariants \`a droite sur la cog\`ebre
oppos\'ee: $F^{op}$ oppos\'ee \`a $F $.

\begin{thm}[Th\'eor\`eme de Dualit\'e .]
Soient $K$ et $E$ deux alg\`ebres associatives de dimensions
finies, avec unit\'es et
$L$ et $F$ les cog\`ebres duales.
Soit une application lin\'eaire $x:L\rightarrow Invd(F)$ et
$y:F\rightarrow Invd(L)$ l'application obtenue par
transposition de l'application $x$. Soient $\pi_{x}$ et
$\pi_{y}$ les mophismes de big\`ebres associ\'es:
$$\pi_{x}:T(L)\rightarrow
U_{x}(\Delta_{L},\epsilon_{L})\subset Invd(T(F)) \ ,et\
\pi_{y}:T(F)\rightarrow V_{y}(\Delta_{F},\epsilon_{F})\subset
Invd(T(L))\ .$$
Alors:\\A) pour tout $w\in T(L) $ et tout $ z\in T(F)$nous
avons:
$$\epsilon_{F} \circ \pi_{x}(w)(z) = \epsilon_{L}\circ
\pi_{y}(\tau (z))(\tau (w)) . 
$$
o\`u $\tau $ est l'anti-automorphisme d'alg\`ebre tensorielle
correspondant au renversement de l'ordre des tenseurs .\\
\\
B) Un \'el\'ement $w\in T(L)$ est une relation dans $U_x $
si et seulement si : $$\pi_{x}(w)(z) =0 ,\ \forall z\in
T(F)\ ,
$$ ou $$\epsilon_{F}\circ \pi_{x}(w)(z) =0 ,\ \forall z\in
T(F)\ , $$ ou $$\epsilon_{L}\circ \pi_{y}(z)(\tau( w)) =0 ,\
\forall z\in T(F) \ .$$
C) Un \'el\'ement $z\in T(F)$ est une relation
dans $V_{y}$ si et seulement si: $$\pi_{y}(z)(w)=0,\ \forall
w\in T(L) \ , $$ ou $$\epsilon_{L}\circ \pi_{y}(z)(w)=0, \
\forall w \in T(L) \ ,$$ ou  $$\epsilon_{F}\circ
\pi_{x}(w)(\tau (z)) =0,\  \forall w\in T(L)\ . $$

\end{thm}

\vspace{0.3cm}
\subsection{  Dualit\'e entre les big\`ebres r\'ealis\'ees
   $U_{x}(\Delta_L , \epsilon_L )$ et
$V_{y}^{op}(\Delta_F ,
\epsilon_F )$.}
\begin{defin}
La dualit\'e canonique entre un op\'erateur invariant
\`a droite ou \`a gauche $a$, agissant sur une cog\`ebre $ C
$ est $<a,c> =\epsilon_C \circ a (c),\  c\in C $ .
\end{defin}
Pour deux cog\`ebres de dimensions finies, quand une
application lin\'eaire
$x$ d'une cog\`ebre
$L$ dans les op\'erateurs invariants \`a droite sur une
cog\`ebre $F$ est donn\'ee , l'application transpos\'ee est
donc du m\^eme type et donc l'on construit les big\`ebres
$U_x $ et $V_y $ ; pour tout $u\in U_x $ , $u$ est contruit
par la donn\'ee d'une classe d'\'equivalence $w\in
T(L)/I_{x}$ et $u=\pi_{x}(w) $; le coproduit $\Delta_L $ est
bien d\'efini par morphisme d'alg\`ebres. De plus dans
\cite{mo3} il est  montr\'e que l'id\'eal $I_{x}$ est aussi
un coid\'eal,  (r\'eunion de coid\'eaux de dimensions finies,
qui sont explicitement d\'ecrits). Il en est de m\^eme pour
la big\`ebre $ V_{y} = \pi_{y}(T(F)) \subset Invd(T(L)) $;
les classes d'\'equivalence donnant les \'el\'ements $v\in
V_{y}$, sont aussi d\'efinies part un id\'eal $I_{y}\subset
T(F) $ qui est aussi un coid\'eal, et  $\tau
(I_{y})=I_{y}^{op}$, est donc, aussi
un coid\'eal pour le m\^eme coproduit
$\Delta_{F}$, ou le coproduit oppos\'e $\Delta_{F}^{op} $ .On
note
$V_{y}^{op}(\Delta_{F},\epsilon_{F})$ la big\`ebre
quotient: $T(F)/I_{y}^{op} $ .\\
a) C'est une big\`ebre r\'ealis\'ee :
$V_{y}^{op}(\Delta_{F},\epsilon_{F}) = \pi_{y} \circ \tau
(T(F)) $.\\
b) Soit un
\'el\'ement
$v^{op}\in T(F)/I_{y}^{op}= V_{y}^{op}$ et un \'el\'ement
$u\in U_{x}$; \\ $<u,I_{y}^{op}>= <\pi_{y}(I_y ),\tau (w)> =0
$, et le couplage $<u,v^{op}> $   entre $ U_x $ et
$V_{y}^{op}$ est bien d\'efini ; il est non d\'eg\'en\'er\'e.
 En effet pour
$u\neq 0 $ il existe $z \in T(F)\ $ tel que
$\epsilon_{F}\circ \pi_{x}(u) (z)\neq 0 $; les \'el\'ements
de
$V^{op}$ \'etant de la forme $\tau (z +I_{y})$  on a:
$<u,\tau (\tau z +I_{y})> \neq 0 $.\\
De m\^eme si $z \neq 0 \in T(F) /I_{y}^{op} $ , $\pi_{y}(\tau
(z)) \neq 0 $, et il existe $u \in U_{x} $ tel que $ <u,z>
\neq 0 $.\\
Par construction dans le cas des big\`ebres $U_x $ et 
$V_{y}^{op}$, l'on a :
\begin{eqnarray}  <u,z_1 .z_2 >& =& 
\sum_{k}<u_{k}^{'},z_1  >.<u_{k}^{''},z_ 2  >
 , \\
 <u_1 .u_2 , z>& =& \sum_{k}<u_2 ,z_{k}^{'}>.<u_1
,z_{k}^{''}> \ . 
\end{eqnarray}
Il suffit de d\'emontrer (2) .D'apr\'es  le th\'eor\`eme de
dualit\'e et les faits que $\epsilon_{U}$, et $\Delta_{V}$ ,
sont des morphismes d'alg\`ebres l'on a :\\  
$<u_1 .u_2 , z> =
\epsilon_U
\circ
\pi_{y}(\tau (z)) (\tau (u_2) .\tau (u_1 ))
$ \\$=$ $\sum_{k}\epsilon_{U} \circ \pi_{y}(\tau
(z_{k}^{'}))(\tau (u_2 )).\epsilon_{U}\circ \pi_{y}( \tau
(z_{k}^{''}))(\tau (u_1))$ ;\\
ce qui donne le r\'esultat en utilisant \`a nouveau la
dualit\'e. \\
Les th\'eor\`emes ci-dessus ainsi que les notations sont 
largement suffisants pour l'exposition de la construction
de l'alg\`ebre enveloppante des champs de vecteurs sur $C^{p}
$ et de son dual, ainsi que pour d\'emontrer et interpr\'eter
l'int\'egration de ceux-ci .\\ Mais pour motiver les
restrictions que nous allons devoir apporter 
\`a la donn\'ee
$x:L\rightarrow Invg,r(F_{a}) $, 
pour avoir une d\'efinition de la notion
de champs de vecteurs associ\'es \`a  un module induit
dans un cadre plus g\'en\'eral, nous poursuivons
l'exposition du probl\`eme de l'\'etude des relations sur
$T(F_{a})$ induites par la donn\'ee d'une application
lin\'eaire
$x$ d'une cog\`ebre $L$ dans $Invg,r(F_{a})$.  
\subsection{Induction de relations sur un module
alg\'ebrique et 
notion de champs de vecteurs.} 
Soit $ x:L
\rightarrow Invg,r(F_{a})$, une application lin\'eaire d'une
cog\`ebre de type fini dans les op\'erateurs invariants \`a
gauche r\'eguliers sur la cog\`ebre approch\'ee $F_a $, et
soient 
$X(l)$ les op\'erateurs invariants
\`a gauche agissant sur $T(F_{a})$, $$X(l)\in
Invg(T(F_{a}))
\subset End(T(F_{a}))\ ,$$ donn\'es par le th\'eor\`eme
2.1 , et $U_x (\Delta_{L},\epsilon_{L})$ la big\`ebre
engendr\'ee par les op\'erateurs $X(l)$ et l'identit\'e dans
$Invg(T(F_{a}))$. On rappelle que ces op\'erateurs peuvent
\^etre consid\'er\'es comme des op\'erateurs invariants \`a
droite $x(l) \in Invd,r( F_{a}(\Delta_{F}^{op})) $ et
respectivement $X(l) \in Invd(T( F_{a}(\Delta_{F}^{op})))$.\\
\\
Cette construction nous a donn\'e donc une repr\'esentation
de l'alg\`ebre $T(L)$ sur
$T(F_{a})$ et le morphime $\pi_{x}:T(L) \rightarrow
Invd(T(F_{a}(\Delta_{F}^{op})))$ nous a donn\'e une
repr\'esentation qui est une action de la big\`ebre libre
$T(L)$ sur l'alg\`ebre $T(F_{a})$: \\
$w\in T(L) \ ,  et\  z_1 , z_2 \in T(F_{a})$ nous avons:
$$\pi_{x}(w)(z_1 .z_2 ) = \sum_{k}\pi_{x}(w_{k}^{'})(z_1
).\pi_{x}(w_{k}^{''})(z_2 )\ ,\ ou\ \Delta_{L}(w)=
\sum_{k}w_{k}^{'}\otimes w_{k}^{''}\ .$$
Par d\'efinition une relation de l'alg\`ebre $U_{x}$ est un
\'el\'ement $w$ de l'alg\`ebre $T(L)$ tel que
$\pi_{x}(w)=0\in U_{x}\subset End(T(F_{a})) \ .$\\
\\
Le probl\`eme est maintenant de munir un module
alg\'ebrique de $F_{a}$ de relations accessibles.\\
Un module alg\'ebrique de $F_{a}$ pour
l'action des op\'erateurs invariants \`a gauche  est obtenu
par un vecteur ligne de la cog\`ebre approch\'ee
$f_{n}^{m}$,
$m$ \'etant fix\'e .\\ Soient
$L_0 $, une cog\`ebre de dimension finie et 
$x_{0}$ une application lin\'eaire $ x_0 :L_0
\rightarrow Invg,r(F_a ).$\\ Par d\'efinition pour $l\in L_0
$ nous avons : 
$$x_{0}(l) (f_{n}^{m}) =\sum_{k}f_{k}^{m}f_{n}^{k}(O(l)) \
.$$ Soit l'application:
$x_{0,t} :l\rightarrow O(l) $;
par transposition on obtient une application :
$y_{0} :F_a \rightarrow K_0$, o\`u $K_0 $ est l'alg\`ebre
duale de la cog\`ebre $L_0 $. \\
\\
La d\'efinition d'une relation dans
$T(F_{a})$  est la suivante.
\begin{defin}

$z\in T(F_{a})$ est une relation si et seulement si

$$\epsilon_{F_{a}}\circ \pi_{x_{0}} (w)(z)= 0\ , \ \forall w
\in T(L_{0})$$
\end{defin}
Mais ce type de relations, n'est pas accessible et
int\'eressant directement ; pour cela il faut restreindre la
donn\'ee de l'application lin\'eaire \\ $ x_0 :L_0
\rightarrow Invg,r(F_a )$ .
\begin{defin}{Module induit. }
On dira qu'une application lin\'eaire \\ $ x_0 :L_0
\rightarrow Invg,r(F_a )$, $L_{0}$ \'etant une cog\`ebre de
dimension finie, induit sur un module alg\'ebrique une
structure d'alg\`ebre si et seulement si, il
existe une suite
 strictement croissante $ (q_{\alpha})_{\alpha} \in N$ telle
que les op\'erateurs $x_{0}(L_{0}) $ laissent les cog\`ebres
$F_{q_{\alpha}}$ invariantes pour tout $\alpha$ , o\`u
$F_{q}$ d\'esigne la cog\`ebre
$(f_{n}^{m}\ ,\ \vert m \vert \leq q \ , \  \vert n \vert
\leq q \ )$ avec $\vert n \vert = \sup_{i\in (1,..,p)}
n_{i}$  .\\

\end{defin}
Dans ce cas les relations dans $ T(F_{a})$ sont accessibles
par le th\'eor\`eme de dualit\'e par l'interm\'ediaire des
relations entre les \'el\'ements d'une cog\`ebre
$F_{q_{\alpha}}$. En particulier les relations entre les
\'el\'ements d'un module alg\'ebrique seront donc
accessibles; d'autre part ce sont les seules qui seront bien
pr\'eserv\'ees par d\'efinition  dans la notion suivante de
big\`ebre enveloppante associ\'ee
\`a un module induit.

\begin{defin}[Big\`ebre enveloppante associ\'ee \`a un 
module induit.] 
Une application lin\'eaire $x_{0}:L_{0}\rightarrow
Invg,r(F_{a}(\Delta_{F}))$ \'etant donn\'ee, satisfaisant
les conditions ci-dessus, soit
$M_{k}$ un module alg\'ebrique de $F_{a}$, fix\'e,  
pour l'action des op\'erateurs invariants \`a gauche , et
l'id\'eal : $I_{y_{0},M_{k} }$, associ\'e par la construction
pr\'ec\'edente dans $T(M_{k}) \subset T(F_a )$ ; on appellera
big\`ebre enveloppante de ce module induit $M_{k,y_{0}}$,
la r\'eunion de toutes les big\`ebres
$U_{x_{\beta}}$, construites \`a partir de la donn\'ee d'une
cog\`ebre $L_{\beta}$ et d'une application 
 $x_{\beta}:L_{\beta} \rightarrow  
Invg,r(F_{a}(\Delta_{F}))$, telles que l'action de la
big\`ebre $U_{x_{\beta}}$ sur $T(F_a )$ laisse l'id\'eal
$I_{y_{0},M_{k}}$ invariant .

\end{defin}
Nous allons pour conclure ce paragraphe th\'eorique 
d\'ecrire des propri\'et\'es suffisantes, sur  la structure
des relations d'un module induit qui permettent de
pr\'eciser les conditions pour construire des applications
$x_{\beta}:L_{\beta}\rightarrow Invg,r(F_{a}) $, telles que
les op\'erateurs $X_{\beta}(l)$  pr\'eservent  les
relations du module induit .

\subsection{Une condition suffisante pour qu'un
  module induit admette des champs de
vecteurs .
}

Soit la donn\'ee d'une application lin\'eaire $x_{0}$ d'une
cog\`ebre $L_{0}$ de dimension finie dans les op\'erateurs
invariants \`a gauche sur la cog\`ebre approch\'ee $
(F_{n}^{m}\ ,\ n,m\in N^{p}) $, satisfaisant aux conditions
de la d\'efinition 2.7 . On d\'esignera cette cog\`ebre
approch\'ee par
$F_{a}(\Delta_{F})$ , et on consid\`ere le module
alg\'ebrique

$M_{o}= ( f_{n}^{o}\ , n \in N^{p} \ ,\  o= (0,..,0)\in
N^{p})$.\\ La donn\'ee $x_{0}: L_{0}\rightarrow Invg,r(F_{a})
$ d\'efinit,
\\a) les op\'erateurs  $ x_{0}(l) (f_{n}^{m})=
\sum_{k}f_{k}^{m}.f_{n}^{k}(O_{r}(l)) \in Invg,r(F_a )$ \\b)
 les op\'erateurs $X(l)\in U_{x_{0}}
(\Delta_{L_{0}},\epsilon_{L_{0}} )$;
\\c)
 l'id\'eal des relations associ\'ees au module induit
: $M_{o,y_{0}}=T(M_{o})/I_{y_{0},M_{o}} $.\\
\\
Dans la cog\`ebre $ L_{0}$  soit le sous espace vectoriel
$L^0
$ ,
$L^0 = ker \epsilon  $ le noyau de la counit\'e .
On suppose de plus: $L_{0} = L^{1}\oplus
L^{0}$ avec $L^{1}$ sous cog\`ebre de $L_{0}$ et que
$x_{0}(l)(f_{o}^{o}) = \epsilon_{L_{0}}(l).f_{o}^{o}\ ,\
l\in L_{0}
$
 de sorte que $f_{o}^{o} $ soit identifi\'e \`a l'unit\'e du
module induit .\\
\\ 
L'approche naturelle de la construction d'applications
$x_{\beta}$ intervenant dans la d\'efinition de l'alg\`ebre
enveloppante
 d'un module induit consiste dans la  construction 
d'applications lin\'eaires 
$$ \beta:L^{0} \rightarrow Hom(M_{0}, T(M_{o})) \ .$$
Soit l'application $x_{\beta}$ :
  
$x_{\beta}:L_{0}\rightarrow Hom(M_{o},T(M_{o}))$ d\'efinie
par:

$$ x_{\beta}(L^1 ) = x_{0}(L^1) , \ et\  \forall l\in L^{0} \
,\  x_{\beta}(l) = \beta(l) \ .
$$
D'apr\'es le Lemme 2.1  de \cite{mo1} il existe une
application  lin\'eaire unique $X:L_{0}\rightarrow
End(T(M_{o}))$  qui satisfasse :\\
a) $$X_{}(l) (1) =\epsilon (l) .1\ ,\ 1 \in T(M_{o}) \ ,$$ 
b) $$X(l) (M_{0}) = x_{\beta}(l)(M_{o})$$
c) $$X(l)( m_1 . m_2 ) =\sum_{k}
X(l_{k}^{'})(m_1).X(l_{k}^{''})(m_2 ) \ ,\ m_1 \ ,\ m_2 \in
T(M_{o}) \ .$$
D'autre part consid\'erons l'espace vectoriel $V$ de
dimension
$p$ :\\
$$V=
(f_{(1,0..0)}^{o},f_{(0,1,0..,0)}^{o},..,f_{(0,..,1)}^{0})$$
Notons ces \'el\'ements $f_{i}\ ,\ i\in (1,2,..,p)$ .
Soit $ T(V)\subset T(M_o) $ et supposons qu'il existe un
sous espace vectoriel de $T(V) \subset T(M_{o})$ que l'on
prendra pour simplifier 
 engendr\'e par les vecteurs de la forme :
$$(f_{1})^{n_{1}}\otimes..\otimes (f_{p})^{n_{p}} =f^{n} $$
 et supposons qu'ils v\'erifient :\\
A) ces \'el\'ements sont lin\'eairement ind\'ependants
modulo 
$I_{y_{o},M_{o}}$;\\
B) c'est un syst\`eme de g\'en\'erateurs pour
$T(M_{o})$ modulo $I_{y_{o},M_{o}}$;\\
C) tout \'el\'ement de cet espace vectoriel :
$$\sum_{n\in N^{p}}a_{n}f^{n}$$
o\`u les coefficients sont presque tous nuls, est dans
$M_{o}$ modulo $I_{y_{o},M_{o}}$ .\\
\\
Sous ces hypoth\`eses A,B,C sur l'id\'eal $I_{y_{0},M_{o}} $
la donn\'ee d'une application : $\beta : L^{0}\rightarrow
Hom(V,T(V)) $ permet de d\'efinir une application lin\'eaire 
$x_{\beta}:L_{0}\rightarrow Invg,r(M_{o})\subset
Invg,r(F_{a})$ telle que les op\'erateurs $X_{\beta}(l) \ ,
\ l\in L_{0}$, alors construits pr\'eservent l'id\'eal
$I_{y_{0},M_{o}}
$ .\\
\\
Remarque : L'hypoth\`ese B) est \'evidemment forte mais peut
\^etre l\'eg\`erement affaiblie.\\
\\
Ce type de construction n'est \'evidemment pas n\'ecessaire
dans le cas des champs de vecteurs sur $C^{p}$, parce que
l'on   dispose de formules combinatoires \'el\'ementaires
qui permettent d'utiliser directement le th\'eor\`eme 2.2
et la d\'efinition 2.6 .\\
\\
 Par contre l'int\'er\^et r\'eside dans la
dualit\'e  canonique qui existe alors entre l'alg\`ebre
enveloppante du module induit et $T(F_{a})$ et donc
aussi
$T(M_{o})$,  ce qui va nous permettre de d\'emontrer des
th\'eor\`emes d'existence et de les interpr\'eter comme
\'etant les  th\'eor\`emes de Cauchy d'int\'egration des
champs de vecteurs. La cog\`ebre approch\'ee jouant  ici un
r\^ole de "cog\`ebre" universelle; en particulier le r\^ole
de son "coproduit" est essentiel et, donne ainsi dans ce cas
pr\'ecis une interpr\'etation et une d\'emonstation purement
alg\'ebrique des th\'eor\`emes d'int\'egration des champs de
vecteurs.\\ 
La raison essentielle pour laquelle  la partie th\'eorique
de cet article a \'et\'e developp\'ee dans une optique non
co-commutative, r\'eside dans les faits que le th\'eor\`eme
de dualit\'e, associ\'e \`a la m\'ethode de r\'ealisation de
big\`ebres
\`a partir de la donn\'ee d'une application lin\'eaire
d'une cog\`ebre $L$ \`a valeurs dans les op\'erateurs
invariants \`a droite (ou \`a gauche ) sur une cog\`ebre
appoch\'ee sont essentiellement des th\'eor\`emes
g\'en\'eraux d'alg\`ebre \cite{mo1}, \cite{mo2}, \cite{mo3} ,
mais qui permettent de donner dans ce cadre une
interpr\'etation et une d\'emonstration des th\'eor\`emes
d'int\'egration des champs de vecteurs sur
$C^p $; la structure de cog\`ebre approch\'ee $F_a $  est
essentielle pour ce qui concerne l'existence, et la
dualit\'e pour ce qui concerne l'interpr\'etation.
\section{Alg\`ebre enveloppante et dual d\'efinis par
les champs de vecteurs sur $C^p $ et leur
int\'egration.}
\subsection{Module induit $M_{o,y_{0}}$ .} 
Soit la cog\`ebre approch\'ee $F_{a}=(f_{n}^{m}\ ,\ m,n \
\in N^{p})$ , et le module alg\'ebrique
$M_{o}=(f_{n}^{o}
\ ,\ o=(0..0),n\ \in N^{p}) $. On veux induire \`a partir de
l'application lin\'eaire $x_{0}:L_{0} \rightarrow Invg,r(F_a
)$ o\`u
$L_{0} $ est la cog\`ebre de Leibnitz de dimension $p+1 $  ,
$l_0,l_1,  ..,l_p
\ ,\
\Delta l_0 =l_0 \otimes l_0 \ ,\ \Delta l_i = l_0 \otimes
l_i +l_i \otimes l_0 \ , \  i\neq 0 $ et par le choix des
op\'erateurs $x_{0}(l)\in Invg,r(F_a )$ la structure
d'alg\`ebre commutative sur le module $M_o $ v\'erifiant les
relations :
$f_{n}^o \otimes f_{m}^o = f_{n+m}^o $.\\
Notations:\\
Soit $n\in N $ on d\'esignera par $i(n) \in N^p $
l'el\'ements
$(0,..0,n,0,..,0 )$ , n \'etant situ\'e , \`a la i\`eme
place.\\
Consid\'erons l'application
$$ x_{0,t} (l_{i})= \sum_{n\in N}  n. E_{i(n)}^{i(n-1)} \ .
$$ Les op\'erateurs invariants \`a gauche r\'eguliers sur
$F_a $, $x_{0}(l_{i})$ transforment :
$$x_{0}(l_i)(f_{n}^{k}) =n_{i} f_{n-i(1)}^{k}$$
$$x_{0}(l_0)= Id  \ .$$
\begin{prop}
L'application $x_{0}:L_{0}\rightarrow Invg,r(F_a )$ ci-dessus
est l'unique application $x_{0}$ d\'efinie sur la cog\`ebre
de Leibnitz qui induisent les relations :
$$f_{n}^{o}\otimes f_{m}^{o}=f_{m+n}^{o} $$ et qui
satisfasse sur les g\'en\'erateurs:
$x_{0}(l_i )f_{o}^{o} =0$ et $x_{0}(l_i ) f_{j(1)}^{o} =
\delta (i,j).f_{o}^{o}$ pour $i >0 $.

\end{prop}
D\'emonstration. 
Soient les op\'erateurs invariants \`a gauche $d_i =X(l_i ) $
donn\'es par le th\'eor\`eme 2.1; ils agissent sur $T(F_a )$
comme d\'erivations et il est facile de montrer qu'ils
commuttent entre eux. Dans notre cas  (th\'eor\`eme 2.2.C) la
cog\`ebre $L_{0}$ \'etant co-commutative, la structure
induite sur
$T(F_a )$ est commutative. Calculons donc: 
$$\epsilon_{F} \circ d_{1}^{n_{1}}\circ
..\circ d_{p}^{n_{p}}(f_{\alpha}^{o} \otimes f_{\beta}^{o} -
f_{\alpha + \beta}^{o} \ )\ . $$

D'une part l'on a :$$\epsilon_{F} \circ d^{n}
f_{\alpha}^{o}= \delta (n,\alpha ) \alpha !     $$

D'autre part l'on a :
$$\Delta_{L}d^{n}= \sum_{p+q=n} n!/p!q!.d^{p}\otimes d^q    
$$

Ce qui montre que l'on a bien les relations attendues sur
le module induit: $M_{o,y_{0}}\subset
T(F_a)/I_{y_{0}}$; il y a pour ce qui concerne $T(F_a )$
d'autres relations int\'eressantes induites par l'actions des
op\'erateurs invariants \`a gauche $d_{i} \ , \ i\in
(1,..,p)$ ; mais ces relations ne seront pas pr\'eserv\'ees
en g\'en\'eral par l'action des champs de vecteurs
associ\'es au sous module $M_{o,y_{0}}$; les relations
pr\'eserv\'ees seront par d\'efinition  et construction en
particulier celles de $M_{o,y_{0}}$ .

\subsection{Sous big\`ebre enveloppante engendr\'ee par
les champs de vecteurs sur $C^{p} .$}
Construction des champs de vecteurs dans le cadre du module 
$M_{o,y_{0}}$.

Pour $i\in (1,2,..,p)$ soient: 
$$A^{i}= \sum_{m\in N^{p}} a_{m}^{i}f_{m}^{o} \in
M_{o}\subset F_{a}\ ;$$
\`a proprement parler ou pour simplifier la lecture on doit
ou l'on peut supposer, que les $a_{m}^{i} $ sont 
 nuls except\'es pour un nombre fini de valeurs $m$ ; on
dira que le champ de vecteurs est analytique , sur un ouvert
$O\subset C^{p}$ si les s\'eries :$$ A^{i}(z)=
\sum_{m\in N^{p}}a^{i}_{m}z_{1}^{m_{1}}..z_{p}^{m_{p}} $$ 
d\'efinissent des fonctions analytiques sur cet ouvert
.\\ Consid\'erons l'application $x_{A}$ de la cog\`ebre
$ L_{0}$ dans les op\'erateurs invariants \`a gauche
r\'eguliers sur
$F_a $ d\'efinie par:  $x_A (l_o) =id$  et
$$x_{A}(l_{i}) (f_{n}^{o} )=n_{i}(\sum_{m\in N^{p}
}a^{i}_{m}.f_{n+m-i(1) }^{o} ) \ .                 
$$
Il est clair que ceci d\'efinit uniquement les
op\'erateurs invariants \`a gauche  $x_{A}(l_i )$ sur $F_a $
et par constuction $X_{A}(l_{i}) $ sur $T(F_{a}) $ :

$$x_{A,t}= \sum_{n\in N^{p} }  \sum_{i\in
(1,..,p)}n_{i}.\sum_{m}a_{m}^{i} E_{n}^{n+m-i(1)}
$$
Soit la big\`ebre
$U_{x_{A}}(\Delta_{L_{0}},\epsilon_{L_{0}})\subset Invg(T(F_a
))$, donn\'ee par le th\'eor\`eme 2.1 dans le cas r\'egulier;
l'on a alors :
\begin{prop}
L'alg\`ebre $U_{x_{A}}$ engendr\'ee dans $End(T(F_a ))$ par
les op\'erateurs invariants \`a gauche r\'eguliers $X_{A}(l_i
)
\in U_{x_{A}}(\Delta_{L_{0}},\epsilon_{L_{0}})$ pr\'eserve
les relations du module induit $M_{o,y_{0}} $ :\\
$$f_{p}^{o}\otimes f_{q}^{o} = f_{p+q}^{o} $$

\end{prop}
D\'emonstration. Calculons : $$ X_{A}(l_i )(f_{p}^{o}\otimes
f_{q}^{o} -f_{p+q}^{o}) \ ;$$ parce que les
op\'erateurs $x_{A}(l_i )$ agissent par construction par
d\'erivation, la premi\`ere partie de l'expression
s' \'ecrit :
$$p_{i}(\sum_{m\in N^{p}}a^{i}_{m}f_{p+m-i(1)}^{o})\otimes
f_{q}^{o} + f_{p}^{o}\otimes q_{i}(\sum_{m\in
N^{p}}a^{i}_{m}f_{q+m-i(1)}^{o}) \ ; 
$$
d'autre part la deuxi\`eme partie s'\'ecrit :

$$( p_{i}+q_{i})(\sum_{m\in N^{p}}
a_{m}^{i}.f_{p+q+m-i(1)}^{o}) $$
Et ces deux expressions coincident bien modulo les
relations du module induit;  elles sont donc
pr\'eserv\'ees . 
\begin{defin}{Champs de vecteurs de composantes
$(A^{i})_{i\in (1,..,p)}$.} C'est par d\'efinition la
d\'erivation
$D_{A}$ obtenue par 
$D_{A}= \sum_{i\in (1,..,p)} X_{A}(l_{i}) $.

\end{defin}
 
\begin{defin}{Big\`ebre enveloppante et son dual associ\'es
aux champs de vecteurs .}
\end{defin}
Consid\'erons la cog\`ebre $L$ somme directe des cog\`ebres
$ L_{A}=L_{0}$, de dimension $p+1$, permettant de
repr\'esenter l'action des champs de vecteurs sur le module
$M_{o,y_{0}}$ par des op\'erateurs invariants \`a gauche
r\'eguliers sur $ T(F_a)$;  on d\'efinit ainsi une
application lin\'eaire:
$$ x=\oplus_{A}\  x_{A}:  L \rightarrow Invg,r 
(F_{a})\ , L=\oplus_{A} L_{A} \ ,\ L_{A}=L_{0}\ ,
$$
et la big\`ebre $
U_{x}(\Delta_{L},\epsilon_{L})\subset
Invg(T(F_{a}))$ admet donc un coproduit et une counit\'e
bien d\'efinis. De plus
$T(F_a ) $ dont les seules relations bien \'etablies sont
les relations polynomiales du module $M_{o,y_{0}}$ , joue
donc le r\^ole de dual pour de la big\`ebre engendr\'ee par
les op\'erateurs $D_{A} \in U_{x}$ .
\begin{thm} 
Soit $\epsilon _{F}$ la counit\'e sur $T(F_{a})$;
pour toute suite d' \'el\'ements 
$D_{A_{1}},D_{A_{2}},..,D_{A_{n}}$ dans la big\`ebre
$U_{x}(\Delta_{L},\epsilon_{L})$ agissant sur $T(F_a )$ l'on
a , en particulier :
$$\epsilon_{F}( D_{A_{1}} \circ D_{ A_{2}}..\circ D_{A_{n}}
(f_{\beta}^{\alpha})) =
\sum_{k_{1}}\sum_{k_{2}}..\sum_{k_{n}}\epsilon_{F}(D_{A_{1}}
(f_{k_{1}}^{\alpha})).\epsilon_{F}(D_{A_{2}}
(f_{k_{2}}^{k_{1}}))..
\epsilon_{F}(D_{A_{n}}(f_{\beta}^{k_{n}} )).$$ Les
sommations
\'etant convergentes :
$$\vert \epsilon_{F}( D_{A_{1}} \circ D_{ A_{2}}..\circ
D_{A_{n}} (f_{\beta}^{\alpha})) \vert \leq  \frac{(deg(\alpha
)+n)!}{deg(\alpha)!}.m(A_{1}).m(A_{2})...m(A_{n}).deg(\beta
)\
\ .
$$
o\`u :
$$\deg(\alpha )= \sum_{i\in (1,..,p)} \alpha_{i}\, \ et \ \ 
m(D_{A})=\sup_{i\in (1,..,p)}\sum_{m}\vert
a_{m}^{i}\vert  .
$$
Si $deg(\beta )> deg(\alpha )+ n+1 $ , alors l'expression est
nulle.

\end{thm}
Ce th\'eor\`eme illustre en particulier l'int\'er\^et 
de la construction de big\`ebres par la donn\'ee d'une
application lin\'eaire d'une cog\`ebre \`a valeurs dans les
op\'erateurs invariants \`a gauche sur une cog\`ebre
approch\'ee , ce qui permet d'interpr\'eter de fa\c{c}on
universelle la loie de composition dans une alg\`ebre
enveloppante en utilisant le coproduit dans le dual qui
s'exprime en terme du coproduit dans la cog\`ebre
approch\'ee.
\\
D\'emonstration.\\
Calculons $\epsilon_{F}(D_{A} (f_{k_{2}}^{k_{1}}))$ par
construction \\
 $D_{A}f_{k_{2}}^{k_{1}}=\sum_{i\in
(1,..,p)}X_{A} (l_{i})f_{k_{2}}^{k_{1}} =\sum_{i\in
(1,..,p)}(k_{2})_{i}.\sum_{m\in
N^{p}}a^{i}_{m}f_{k_{2}+m-i(1)}^{k_{1}} \ .
$\\
\\
Remarquons que $\epsilon_{F}(
(k_{2})_{i}.a^{i}_{m}f_{k_{2}+m-i(1)}^{k_{1}} )
$ est non nul seulement si \\
\\
le degr\'e $deg (k_{2}) \leq
deg(k_{1})+1 $  o\`u   \   $deg(k) = \sum_{i\in (1,..,p)}
k_{i}$, $k\in N^{p}$.
\\
Par simplicit\'e introduisons la quantit\'e
suivante, associ\'ee \`a chaque champ de vecteurs $D_{A}$:
$m(D_{A}) = \sup_{i\in (1,..,p)}\sum_{m}\vert
a_{m}^{i}\vert $ .\\
Alors on a les majorations suivantes pour 
$\epsilon_{F}(D_{A} (f_{k_{2}}^{k_{1}}))$ :
 $$\vert \epsilon_{F}(D_{A} (f_{k_{2}}^{k_{1}}))\vert \leq
deg(k_{2}).m(A)
\leq (deg (k_{1})+1) .m(A) .$$ 
Ainsi on obtient:
$$\epsilon_{F}( D_{A_{1}} \circ D_{ A_{2}}..\circ D_{A_{n}}
(f_{\beta }^{\alpha })) < \frac{(deg(\alpha
)+n)!}{deg(\alpha )!}.m(A_{1}).m(A_{2})...m(A_{n}).deg(\beta
)\
\ .
$$
L'expression est nulle si $deg(\beta) $ est
sup\'erieure \`a $deg(\alpha )+n+1$.

\begin{thm}
Soient deux champs de vecteurs analytiques $ A\ ,\ B$ et
$D_{A}\ ,\ D_{B}$ les d\'erivations associ\'ees (def.3.1) ; 
consid\'erons dans l'alg\`ebre enveloppante \'etendue et
$t\in C$ :
$e^{t. D_{A}} = \sum_{n\in N}
\frac{t^{n}}{n!}D_{A}^{n} $. Alors:\\
\\
a) $\epsilon_{F}(e^{t.D_{A}}(f_{\beta}^{o}))$ existe, 
$t\in C $  pour $ \vert t \vert < 1/m(A) \ .$ \\
\\
b)  Soient les champs de vecteurs $A$ et $B$ et les nombres 
$m(A),m(B)$ alors:
$$\epsilon_{F} (e^{t_{2}.D_{B}} \circ e^{t_{1}.D_{A}}
(f_{\beta}^{o}))    
  $$
existe pour : $\vert t_{1}\vert. m(A) +\vert
t_{2}
\vert . m(B) <1 $

\end{thm}

D\'emonstration.\\
Le terme g\'en\'erique de la s\'erie $
\frac{t^{n}}{n!}\epsilon_{F}(D_{A}^{n}(f_{\beta}^{o}))  $
est donc born\'e d'apr\'es le th.3.1 par $(\vert
t\vert)^{n}.m(A)^{n} $ et donc pour $\vert t \vert < 1/m(A) $
la s\'erie est convergente.
De m\^eme $\epsilon_{F} (e^{t_{2}.D_{B}} \circ
e^{t_{1}.D_{A}} (f_{\beta}^{o}))    
  $ est obtenu par la somme double sur $n\in N$ et $m \in N$
dont le terme g\'en\'erique est d'apr\'es le th.3.1 
major\'e par :
$$
\frac{(t_{2}.m(B))^{n}}{n!}.\frac{(t_{1}.m(A))^{m}}{m!}(n+m)!\
,$$

et cette somme double est \'egale \`a :

$$ \sum_{p\in N} (t_{2}.m(B) + t_{1}m(A))^{p} .$$
Ce qui est convergent sous l'hypoth\`ese du thm.3.2  ;
les propri\'et\'es de diff\'erentiabilit\'e en $t_{1}$, et
$t_{2} $ deviennent alors \'evidentes sous les m\^emes
hypoth\`eses.\\ 
{\bf Remarque} \\Il faut noter 
la  tr\'es faible d\'ependance en $deg(\beta )$ intervenant
dans le th\'eor\`eme 3.1 .

\subsection{Interpr\'etation : int\'egration des champs de
vecteurs  sur $C^{p}$ .}

L'interpr\'etation dans le cadre de big\`ebres en dualit\'e
est bas\'ee dans le cas de big\`ebres enveloppantes
construites \`a partir de cog\`ebres de type Leibnitz sur
les faits que l'on obtient des repr\'esentations de
$T(F_a)$ et en particulier du module induit $M_{o,y_{0}}$
\`a partir de coid\'eaux \`a gauche minimaux de
l'alg\`ebre enveloppante que l'on obtient dans notre cas
par l'exponentielle d'une d\'erivation :\\ $\Delta_{L}(
e^{t.D_{A}}) =e^{t.D_{A}}\otimes e^{t.D_{A}}$ .\\ De plus de
mani\`ere g\'en\'erale,
$\epsilon_{F}$ est un morphisme de l'alg\`ebre tensorielle
$T(F_{a})$ dans $C$ . 
Soit $f_{n}^{o}\in M_{o,y_{0}}$ ; en rappelant que
$i(1) $ repr\'esente l'\'el\'ement dans $N^{p} $ :
$(0,.0,1,0..0)$ la valeur $1$ \'etant situ\'ee \`a la i\`eme
place et en notant  $ f_{i}^{0} = f_{i(1)}^{0}\in
M_{o,y_{0}}$  \  on a :
 $ f_{n}^{o} =
(f_{1})^{n_{1}}.(f_{2})^{n_{2}}..(f_{p})^{n_{p}}.$\\
\\
En utilisant les propri\'et\'es mentionn\'ees 
on obtient:\\

$\epsilon_{F}(e^{t.D_{A}}
(f_{n}^{o}))=\epsilon_{F}(e^{t.D_{A}} (f_{1}^{o}))^{n_{1}}.
 \epsilon_{F}(e^{t.D_{A}}
(f_{2}^{o}))^{n_{2}}...\epsilon_{F}(e^{t.D_{A}}
(f_{p}^{o}))^{n_{p}}.$\\
Notons $$y_{i}(t) =
\epsilon_{F}(e^{t.D_{A}}(f_{i}^{o})) \ ,\ pour \ i\in
(1,..,p) .
$$

$\frac{d}{dt}y_{i}(t) =\epsilon_{F}(e^{t.D_{A}}\circ
D_{A}(f_{i}^{o})) =
\epsilon_{F}(e^{t.D_{A}}(\sum_{m}a^{i}_{m}f_{m}^{o}))$ .
Et d'apr\'es les remarques ci-dessus, ceci nous donne: 
 $$\frac{d}{dt}y_{i}(t) = A^{i}(y_{1}(t),y_{2}(t)..,y_{p}(t))
\ ,\ \forall i\in (1,..,p) .$$
\hspace{2.5cm}    $y_{i}(0) =0 \ \  \forall i\in
(1,..p)
$ .\\
Plus g\'en\'eralement soit $(x^{1},..,x^{p}) \in C^{p}$
, en consid\'erant les op\'erateurs invariants \`a gauche
\'el\'ementaires, qui ont donn\'e les d\'erivations $d_{i}$
commuttant entre elles, on construit $
e^{x^{i}.d_{i}}$ et donc l'on construit l'homomorphisme
d'alg\`ebre $
 e^{x^{i}.d_{i}}\circ e^{t.D_{A}}$.
Pour $x$ dans un ouvert appropri\'e 
l'on d\'enote :
$$y_{i,x}(t)= \epsilon_{F}( e^{x^{i}.d_{i}}\circ e^{t.D_{A}}
(f_{i}^{o}))=\epsilon_{F}(e^{t.D_{A_{x}}}\circ
e^{x^{i}.d_{i}} (f_{i}^{o}))$$
o\`u $D_{A_{x}}=e^{x^{i}.d_{i}}\circ D_{A} \circ
e^{-x^{i}.d_{i}}$, et en choisissant $t$ v\'erifiant $\vert
t\vert. m(A_{x}) <1$ pour tout $x$ dans l'ouvert
consid\'er\'e , l'on obtient le th\'eor\`eme suivant :\\
\\
\\

\begin{thm}
Pour tout champ de vecteurs analytiques $A$  pour tout $x\in
O^{p}$ un ouvert \`a adh\'erence compacte de $C^{p}$, il
existe $0\in O\subset C$ un voisinage de l'origine dans
$C$ tel que pour     
$t\in O $, et $x\in O^{p}$ l'on a : $m(A_{x}).\vert t\vert<1
$ et les fonctions :
$$y_{i,x}(t)= \epsilon_{F}(e^{x^{i}.d_{i}} \circ
e^{t.D_{A}}(f_{i}^{o}))$$
v\'erifient: 
$$\frac{d}{dt}y_{i,x}(t) =
A^{i}(y_{1,x}(t),y_{2,x}(t),..,y_{p,x}(t))\ ,
 $$
et 
$$ y_{i,x}(0) =x_{i} \ ,\ \forall i\in (1,..,p) .$$

\end{thm}
Rappelons que pour tout champ de vecteurs
analytiques sur 
$C^{p}$, la quantit\'e $m(A)$ est donn\'ee par ;
$$m(A) =\sup_{i\in (1,2,..,p)} \sum_{m\in N^{p}} \vert
a_{m}^{i}
\vert \ ,$$
o\`u la i-\`eme composante du champ est donn\'ee pour
$z=(z_{1},..,z_{p}) \in C^{p}$ par :
$$A^{i}(z) = \sum_{m\in N^{p}}
a_{m}^{i}.z_{1}^{m_{1}}.z_{2}^{m_{2}}...z_{p}^{m_{p}}
\ .$$
 Remarques :\\
A) Les estimations ci-dessus  sont des estimations
g\'en\'erales valables pour tous les champs de vecteurs.
Nous avons d\'ej\`a signal\'e que dans ce cas nous avions peu
d'informations sur le dual $T(F_{a})$ \`a l'exception de
celles concernant le module induit . Mais s'il s'agit
d'\'etudier un champ de vecteurs particulier la seule
donn\'ee  de l'application lin\'eaire $x_{A}$ de la
cog\`ebre de Leibnitz $L_{A}$ de dimension $p+1$ ou $1+1$ \`a
valeurs dans les op\'erateurs invariants \`a gauche sur la
cog\`ebre approch\'ee $F_{a}$, induit des relations sur
$T(F_{a})$ dont l'\'etude th\'eorique est plus
abordable.\\
\\
B) Remarquons que lorsqu'il existe un sous espace de 
 dimension fini du module $M_o$ invariant par la
d\'erivation $D_{A}$ d\'efinie par le champ de vecteurs,
qu'alors l'int\'egration sur la sous alg\`ebre engendr\'ee
dans $M_o$ par ce sous espace, est ramen\'ee \`a une
exponentielle d'un matrice de dimension finie et donc ne
pr\'esente pas de singularit\'es
\`a temps fini. Cette remarque met en \'evidence les
diff\'erences et les liens dans le cadre de la dualit\'e
expos\'ee entre les alg\`ebres de Lie de dimensions finies ,
et l'alg\`ebre de Lie des champs de vecteurs que l'on
peut inclure dans un m\^eme cadre en introduisant les
op\'erateurs invariants \`a gauche (o\`u \`a droite ) sur
une cog\`ebre approch\'ee et leurs actions sur $T(F_{a})$
d\'efinies par des applications convenables d'une cog\`ebre
$L$ dans $Invg(F_{a})$ . La possibilit\'e de l'existence
de singularit\'es dans l'int\'egration d'un champ de
vecteurs   est intrinsecte \`a la non existence dans le
module d'un sous espace vectoriel de dimension
fini, invariant par l'action de ce champ de vecteurs, comme
le montrent des exemples \'el\'ementaires en dimension
$1$ .

\end{document}